\documentclass[psamsfonts]{amsart}
\usepackage{amssymb,latexsym,cite,enumerate,}
\usepackage[mathscr]{eucal}
\usepackage[cp1251]{inputenc}
\usepackage[english]{babel}

\newtheorem{theorem}{Theorem}

\newtheorem{lemma}{Lemma}
\theoremstyle{definition}

\newtheorem{corollary}{Corollary}
\theoremstyle{remark}

\renewcommand{\P}{\mathbf{P}\,}
\newcommand{\E}{\mathbf{E}\,}
\newcommand{\N}{\mathbb{N}}
\newcommand{\Z}{\mathbb{Z}}
\newcommand{\R}{\mathbb{R}^1}

\allowdisplaybreaks
\begin{document}
\title{On distribution of zeros of random polynomials in complex plane}
\author{Ildar Ibragimov and Dmitry Zaporozhets}
\begin{abstract}

 Let $G_n(z)=\xi_0+\xi_1z+\dots +\xi_n z^n$ be a random polynomial with i.i.d. coefficients (real or complex). We show that the arguments of the roots of $G_n(z)$ are uniformly distributed in $[0,2\pi]$ asymptotically as $n\to\infty$.  We also prove that the condition $\E\ln(1+|\xi_0|)<\infty$ is necessary and sufficient for the roots  to asymptotically concentrate near the unit circumference.
\keywords{roots of random polynomial, roots concentration, random analytic function}

%\emph{Key words and concepts:} roots of a random polynomial, roots concentration, random series.

\end{abstract}

\thanks{Partially supported by RFBR (08-01-00692, 10-01-00242), RFBR-DFG (09-0191331), NSh-4472.2010.1, and CRC 701 ``Spectral Structures and Topological Methods in Mathematics''}
 \maketitle

\section{Inroduction: problem and results}
\label{sec:1}

Let $\{\xi_k\}_{k=0}^\infty$ be a sequence of independent identically distributed real- or complex-valued random variables.  It is always supposed that $\P(\xi_0=0)<1$.

Consider the sequence of random polynomials
$$
G_n(z)=\xi_0+\xi_1z+\dots+\xi_{n-1} z^{n-1} +\xi_n z^n\;.
$$
By $z_{1n},\dots,z_{nn}$ denote the zeros of $G_n$. It is not hard to show (see \cite{jH56}) that there exist an indexing of the zeros such that for each $k=1,\dots,n$ the $k$-th zero $z_{kn}$ is a one-valued random variable. For any measurable subset $A$ of $\mathbb{C}$ put $N_n(A)=\#\{z_{kn}\,:\,z_{kn}\in A\}$. Then $N_n(A)/n$ is a probability measure on the plane (the empirical distribution of the zeros of $G_n$). For any $a,b$ such that $0\leqslant a<b\leqslant\infty$ put $R_n (a,b)=N_n(\{z\,:\,a\leqslant|z|\leqslant b\})$  and for any $\alpha,\beta$ such that 0 $\leqslant\alpha<\beta\leqslant2\pi$ put $S_n(\alpha,\beta)=N_n(\{z\,:\,\alpha\leqslant\arg z\leqslant\beta\})$. Thus $R_n/n$ and $S_n/n$ define the empiric distributions of $|z_{kn}|$ and $\arg z_{kn}$.

In this paper we study the limit distributions of $N_n,R_n,S_n$ as $n\to\infty$.

The question of the distribution of the complex roots of $G_n$ have been originated by Hammersley in \cite{jH56}. The asymptotic study of $R_n ,S_n$ has been initiated by Shparo and Shur in \cite{SS62}. To describe their results let us introduce the function
$$
f(t)=\left[\underbrace{\log^+\log^+\dots\log^+t}_{m+1}\right]^{1+\varepsilon}\prod_{i=1}^m\underbrace{\log^+\log^+\dots\log^+t}_{i}\;,
$$
where $\log^+s=\max(1,\log s)$. We assume that  $\varepsilon>0, m\in\mathbb{Z^+}$ and $ f(t)=(\log^+t)^{1+\varepsilon}$ for $m=0$.

Shparo and Shur have proved in \cite{SS62} that if
$$
\E f(|\xi_0|)<\infty
$$
for some $\varepsilon>0,m\in\mathbb{Z^+}$, then for any $\delta\in(0,1)$ and $\alpha,\beta$ such that 0 $\leqslant\alpha<\beta\leqslant2\pi$
$$
\frac1nR_n(1-\delta,1+\delta)\overset{\P}{\longrightarrow}1\;,\quad n\to\infty\;,
$$
$$
 \frac1n S_n(\alpha,\beta)\overset{\P}{\longrightarrow}\frac{\beta-\alpha}{2\pi}\;,\quad n\to\infty\;.
$$
The first relation means that under quite weak constraints imposed on the coefficients of a random polynomial, almost all its roots ``concentrate uniformly'' near the unit circumference with high probability; the second relation means that the arguments of the roots are asymptotically uniformly distributed.

Later Shepp and Vanderbei \cite{SV95} and Ibragimov and Zeitouni \cite{IZ95} under additional conditions imposed on the coefficients of $G_n$ got more precise asymptotic formulas for $R_n$.

What kind of further results could be expected? First let us note that if, e.g.,  $\E|\xi_0|<\infty$, then for $|z|<1$
$$
G_n(z)\to G(z)=\sum_{k=0}^\infty\xi_kz^k
$$
as $n\to\infty$ a.s. The function $G(z)$ is analytical inside the unit disk $\{|z|<1\}$. Therefore for any $\delta>0$ it has only a finite number of zeros in the disk  $\{|z|<1-\delta\}$. At the other hand, the average number of zeros in the domain $|z|>1/(1-\delta)$ is the same (it could be shown if we consider the random polynomial $G(1/z)$). Thus one could expect that under sufficiently weak constraints imposed on the coefficients of a random polynomial the zeros concentrate near the unit circle $\Gamma=\{z\,:\,|z|=1\}$ and a measure $R_n/n$ converges to the delta measure at the point one. We may expect also from the consideration of symmetry that the arguments $\arg z_{kn}$ are asymptotically uniformly distributed. Below we give the conditions for these hypotheses to hold. We shall prove the following three theorems about the behavior of $N_n/n,R_n/n,S_n/n$.

For the sake of simplicity, we assume that $\P\{\xi_0=0\}=0$.  To treat the general case it is enough to study in the same way the behavior of the roots on the sets $\{\theta'_n =k,\theta''_n=l\}$, where
$$
\theta'_n=\max\{i=0,\dots,n\mid\xi_{i}\ne0\},\qquad
\theta''_n=\min\{j=0,\dots,n\mid\xi_{j}\ne0\}\;.
$$

\begin{theorem}\label{1612}
The sequence of the empirical distributions $R_n/n$ converges to the delta measure at the point one almost surely if and only if
\begin{equation}\label{1427}
\E\log(1+|\xi_0|)<\infty\;.
\end{equation}
\end{theorem}
In other words, (\ref{1427}) is necessary and sufficient condition for
\begin{equation}\label{2149}
\P\left\{\frac1nR_n(1-\delta,1+\delta)\underset{n\to\infty}{\longrightarrow}1\right\}=1
\end{equation}
hold for any $\delta>0$.

We shall also prove that if (\ref{1427}) does not hold then no limit distribution for $\{z_{nk}\}$ exist.

\begin{theorem}\label{1233}
Suppose the condition (\ref{1427}) holds. Then the empirical distribution $N_n/n$ almost surely converges to the probability measure $N(\cdot)=\mu(\cdot\cap\Gamma)/(2\pi)$, where $\Gamma=\{z\,:\,|z|=1\}$ and $\mu$ is the Lebesgue measure.
\end{theorem}

\begin{theorem}\label{1535}
The empirical distribution $S_n/n$ almost surely converges to the uniform distribution, i.e.,
$$
\P\left\{\frac1nS_n(\alpha,\beta)\underset{n\to\infty}{\longrightarrow}\frac{\beta-\alpha}{2\pi}\right\}=1
$$
for any $\alpha,\beta$ such that 0 $\leqslant\alpha<\beta\leqslant2\pi$.
\end{theorem}

Let us remark here that Theorem~\ref{1535} does not require any additional conditions on the sequence $\{\xi_k\}$.

The next result is of crucial importance in the proof of Theorem~\ref{1535}.
\begin{theorem}\label{2101}
Let $\{\eta_k\}_{k=0}^\infty$ be a sequence of  independent identically distributed real-valued random variables. Put $g_n(x)=\sum_{k=0}^n\eta_kx^k$ and by $M_n$ denote the number of real roots of the polynomial $g_n(x)$. Then
$$
\P\left\{\frac{M_n}{n}\underset{n\to\infty}{\longrightarrow}0\right\}=1,\quad \E M_n=o(n),\quad n\to\infty\;.
$$
\end{theorem}

Theorem~\ref{2101} is also of independent interest. In a number of papers it was shown that under weak conditions on the distribution of $\eta_0$ one has $\E M_n\sim c\times\log n, n\to\infty$ (see  \cite{IM71}, \cite{IM71b}, \cite{IM71c}, \cite{mK43}, \cite{LS68}, \cite{LS68b}). L.~Shepp proposed the following conjecture: for any distribution of $\eta_0$ there exist positive numbers $c_1,c_2$ such that $\E M_n\geqslant c_1\times\log n$ and $\E M_n\leqslant c_2\times\log n$ for all $n$. The first statement was disproved in \cite{dZ05}, \cite{ZN08}. There was constructed a random polynomial $g_n(x)$ with $\E M_n<1+\varepsilon$. It is still unknown if the second statement is true. However, Theorem~\ref{2101} shows that an arbitrary random polynomial can not have too much real roots (see also \cite{SF10}).

In fact, in the proof of Theorem~\ref{1535} we shall use a slightly generalized version of Theorem~\ref{2101}:
\begin{theorem}\label{1713}
For some integer $r$ consider a set of $r$ non-degenerate probability distributions. Let   $\{\eta_k\}_{k=0}^\infty$ be a sequence of  independent real-valued random variables with distributions from this set. As above, put $g_n(x)=\sum_{k=0}^n\eta_kx^k$ and by $M_n$ denote the number of real roots of the polynomial $g_n(x)$. Then
\begin{equation}\label{1611}
\P\left\{\frac{M_n}{n}\underset{n\to\infty}{\longrightarrow}0\right\}=1,\quad \E M_n=o(n), n\to\infty\;.
\end{equation}
\end{theorem}

\section{Proof of Theorem~\ref{1612}}
Let us establish the sufficiency of  (\ref{1427}).  Let it hold and fix $\delta\in(0,1)$. Prove that the radius of convergence of the series
\begin{equation}\label{1801}
G(z)=\sum_{k=0}^\infty\xi_kz^k
\end{equation}
is equal to one with probability one.

Consider $\rho>0$ such that $\P\{|\xi_0|>\rho\}>0$. Using the Borel-Cantelli lemma we obtain that with probability one the sequence $\{\xi_k\}$ contains infinitely many $\xi_k$ such that  $|\xi_k|>\rho$. Therefore the radius  of convergence of the series (\ref{1801}) does not exceed 1 almost surely.

On the other hand, for any non-negative random variable $\zeta$
\begin{equation}\label{1215}
\sum_{k=1}^\infty\P(\eta\geqslant k)\leqslant\E\zeta\leqslant 1+\sum_{k=1}^\infty\P(\zeta\geqslant k)\;.
\end{equation}
Therefore, it follows from (\ref{1427}) that
$$
\sum_{k=1}^\infty\P(|\xi_k|\geqslant e^{\gamma k})<\infty
$$
for any positive constant $\gamma$. It follows from the  Borel-Cantelli lemma that with probability one $|\xi_k|<e^{\gamma k}$ for all sufficiently large $k$. Thus, according to the Cauchy-Hadamard  formula (see, e.g., \cite{aM83}), the radius of convergence of the series (\ref{1801}) is at least 1 almost
surely.

Hence with probability one $G(z)$ is an analytical function inside the unit ball $\{|z|<1\}$. Therefore if $0\leqslant a<b<1$, then $R(a,b)<\infty$, where $R(a,b)$ denotes the number of the zeros of $G$ inside the domain $\{z\,:\,a\leqslant |z|\leqslant b\}$. It follows from the Hurwitz theorem (see, e.g., \cite{aM83}) that with probability one $R_n\left(0,1-\delta\right)\leqslant R\left(0,1-\delta/2\right)$ for all sufficiently large $n$. This implies
$$
\P\left\{\frac 1nR_n(0,1-\delta)\underset{n\to\infty}\longrightarrow0\right\}=1\;.
$$
In order to conclude the proof of (\ref{2149}) it remains to show that
$$
\P\left\{\frac 1nR_n(1+\delta,\infty)\underset{n\to\infty}\longrightarrow0\right\}=1\;.
$$
In other words, we need to prove that $\P\{A\}=0$, where $A$ denotes the event that there exists $\varepsilon>0$ such that
$$
R_n\left(1+\delta,\infty\right)\geqslant\varepsilon n
$$
holds for infinitely many values $n$.

By $B$ denote the event that $G(z)$ is an analytical function  inside the unit disk $\{|z|<1\}$. For $m\in\N$ put
$$
\zeta_m=\sup_{k\in\Z^+}|\xi_k e^{-k/m} |.
$$
By $C_m$ denote the event that $\zeta_m<\infty$. It was shown above that $\P\{B\}=\P\{C_m\}=1$ for $m\in\N$. Therefore, to get $\P\{A\}=0$, it is sufficient to show that $\P\{ABC_m\}=0$ for some $m$.

Let us fix $m$. The exact value of it will be chosen later. Suppose the event $ABC_m$ occurred. Index the roots of the polynomial $G_n(z)$ according to the order of magnitude of their absolute values:
$$
|z_1|\leqslant|z_2|\leqslant\dots\leqslant|z_n|.
$$
Fix an arbitrary number  $C>1$ (an exact value will be chosen later). Consider indices $i,j$ such that
\begin{align*}
|z_i|<1-\delta/C\;,&\quad |z_{i+1}|\geqslant 1-\delta/C\;,
\\
|z_j|\leqslant 1+\delta\;,&\quad |z_{j+1}|>1+\delta.
\end{align*}
If $|z_1|\geqslant 1-\delta/C$, then $i=0$; if $|z_{n}|\leqslant 1+\delta$ then $j=n$.

It is easily shown that if
$$
|z|<\min\left(1,\frac{|\xi_0|}{n\times\max_{k=1,\dots,n}|\xi_k|}\right)\;,
$$
then
$$
|\xi_0|>|\xi_1z|+|\xi_2z^2|+\dots+|\xi_nz^n|\;.
$$
Therefore such $z$ can not be a zero of the polynomial $G_n$. Taking into account that the event $C_m$ occurred, we obtain a lower bound for the absolute values of the zeros for all sufficiently large $n$:
$$
|z_1|\geqslant\min\left(1,\frac{|\xi_0|}{n\times\max_{k=1,\dots,n}|\xi_k|}\right)\geqslant\frac{|\xi_0|}{n\zeta_m e^{n/m}}\geqslant|\xi_0|\zeta_m^{-1}e^{-2n/m}\;.
$$
Therefore for any integer $l$ satisfying $j+1\leqslant l\leqslant n$ and all sufficiently large $n$
\begin{multline*}
|z_1\dots z_l|=|z_1\dots z_i||z_{i+1}\dots z_{j}||z_{j+1}\dots
z_l|\\\geqslant |\xi_0|^i\zeta_m^{-i}e^{-2ni/m}\left(
1-\frac\delta C\right)^{j-i}(1+\delta)^{l-j}\;.
\end{multline*}
Since $A$ occurred, $n-j\geqslant n\varepsilon$ for infinitely many values of n. Therefore if $l$ satisfies $n-\sqrt{n}\leqslant l\leqslant n$, then the inequalities $j+1\leqslant l\leqslant n$ and $l-j\geqslant n\varepsilon/2$ hold for infinitely many values of n. According to the Hurwitz theorem, $i\leqslant R_n(0,1-\delta/C)\leqslant R(0,1-\delta/(2C))$ for all
sufficiently large $n$. Therefore for infinitely many values of $n$ 
$$
|z_1\dots z_l|\geqslant\left(\frac{|\xi_0|}{\zeta_m}\right)^{R(0,1-\delta/(2C))}e^{-2nR(0,1-\delta/(2C))/m}\left( 1-\frac\delta C\right)^{n}(1+\delta)^{n\varepsilon/2}\;.
$$

Choose now $C$ large enough to yield
$$
\left(1-\frac\delta C\right)(1+\delta)^{\varepsilon/2}>1\;.
$$
Furthermore, holding $C$ constant choose $m$ such that
$$
b=e^{-2R(0,1-\delta/(2C))/m}\left(1-\frac\delta C\right)(1+\delta)^{\frac\varepsilon2}>1\;.
$$
Since
$$
\left(\frac{|\xi_0|}{\zeta_m}\right)^{R(0,1-\delta/(2C))/n}\underset{n\to\infty}\longrightarrow1\;,
$$
there exists a random variable $a>1$ such that for infinitely many values of n
$$
|z_1\dots z_l|\geqslant\left(\frac{|\xi_0|}{\zeta_m}\right)^{R(0,1-\delta/(2C))}b^n=\left(b\left(\frac{|\xi_0|}{\zeta_m}\right)^{R(0,1-\delta/(2C))/n}\right)^n\geqslant a^n\;.
$$
On the other hand, it follows from $n-\sqrt{n}\leqslant l$ and Vi\'ete's formula that
$$
|z_{l+1}\dots z_n|\geqslant{n \choose n-\sqrt{n}}^{-1}|\sum_{i_1<\dots<i_{n-l}}z_{i_1}\dots
z_{i_{n-l}}|={n \choose n-\sqrt{n}}^{-1}\frac{|\xi_l|}{|\xi_n|}.
$$
We combine these two inequalities to obtain for infinitely many values of n
\begin{multline*}
\frac{|\xi_0|}{|\xi_n|}=|z_{1}\dots z_n|\geqslant a^n{n \choose n-\sqrt{n}}^{-1}\frac{|\xi_l|}{|\xi_n|}\\
\geqslant c_1a^n\frac{(\sqrt n)^{\sqrt n+\frac12}(n-\sqrt n)^{n-\sqrt n+\frac12}}{n^{n+\frac12}}\frac{|\xi_l|}{|\xi_n|}\geqslant c_2a^n(\sqrt n)^{-\sqrt n}\left(1-\frac1{\sqrt n}\right)^{n}\frac{|\xi_l|}{|\xi_n|}\\
\geqslant c_3 \exp\left(n\log a-\frac{\sqrt n\log n}2-\sqrt n\right)\frac{|\xi_l|}{|\xi_n|}\geqslant e^{\alpha n}\frac{|\xi_l|}{|\xi_n|}\;,
\end{multline*}
where $\alpha$ is a positive random variable. Multiplying left and
right parts by $|\xi_n|$, we get
$$
ABC_m\subset\bigcup_{i=1}^{\infty}D_i\;,
$$
where $D_i$ denotes the event that $ |\xi_0|>e^{n/i}\max_{n-\sqrt n\leqslant l\leqslant n}|\xi_l|$ for infinitely many values of $n$.

To complete the proof it is sufficient to show that $\P\{D_i\}=0$ for all $i\in\N$. Having in mind to apply the Borel-Cantelli lemma, let us introduce the following events:
$$
H_{in}=\left\{|\xi_0|>e^{n/i}\max_{n-\sqrt n\leqslant l\leqslant n}|\xi_l|\right\}\;.
$$
Considering $\theta>0$ such that $\P\{|\xi_0|\leqslant\theta\}=F(\theta)<1$, we have
$$
H_{in}\subset\left\{|\xi_0|>\theta e^{n/i}\right\}\cup\left\{\max_{n-\sqrt n\leqslant l\leqslant n}|\xi_l|\leqslant\theta\right\},
$$
consequently,
$$
\sum_{n=1}^\infty\P\{H_{in}\}\leqslant\sum_{n=1}^\infty\P\{|\xi_0|>\theta e^{n/i}\}+\sum_{n=1}^\infty(F(\theta))^{\sqrt n}<\infty
$$
and, according to the Borel-Cantelli lemma, $\P\{D_i\}=0$.

We  prove the implication $(\ref{2149})\Rightarrow(\ref{1427})$ arguing by contradiction. Suppose (\ref{1427}) does not hold, i.e.,
$$
\E\log(1+|\xi_o|)=\infty\;.
$$
It follows from (\ref{1215}) that
\begin{equation}\label{1037}
\sum_{n=1}^\infty\P(|\xi_n|\geqslant e^{\gamma n})=\infty
\end{equation}
for an arbitrary positive $\gamma$. For $k\in\mathbb{N}$ introduce an event $F_k$ that $|\xi_n|\geqslant e^{kn}$ holds for infinitely many values of $n$. It follows from (\ref{1037}) and the Borel-Cantelli  lemma that $\P\{F_k\}=1$ and, consequently, $\P\{\cap_{k=1}^\infty F_k\}=1$. This yields
$$
\P\left\{\limsup_{n\to\infty}|\xi_n|^{1/n}=\infty\right\}=1\;.
$$
Therefore with probability one for infinitely many values of $n$
$$
|\xi_n|^{1/n}>\max_{i=0,\dots,n-1}|\xi_i|^{1/i},\quad|\xi_n|^{1/n}>\frac3\varepsilon,\quad |\xi_0|<2^{n-1},
$$
where $\varepsilon>0$ is an arbitrary fixed value. Let us hold one
of those $n$. Suppose $|z|\geqslant\varepsilon$. Then
\begin{multline*}
|\xi_0+\xi_1z+\dots+\xi_{n-1}z^{n-1}|\\
\leqslant2^{n-1}+|\xi_nz^n|^{1/n}+|\xi_nz^n|^{2/n}+\dots+|\xi_nz^n|^{(n-1)/n}\\
=\frac{2^{n}}{2}-1+\frac{|\xi_nz^n|-1}{|\xi_n^{1/n}z|-1}\leqslant\frac{|\xi_n^{1/n}z|^n}{2}-1+\frac{|\xi_nz^n|-1}{(3/\varepsilon)\times\varepsilon-1}<|\xi_nz^n|\;.
\end{multline*}
Thus with probability one for infinite number of values of $n$ all the roots of the polynomial  $G_n$ are located inside the circle  $\{z\,:\,|z|=\varepsilon\}$, where $\varepsilon$ is an arbitrary positive constant. This means that (\ref{2149}) does not hold for any $\delta\in(0,1)$.

\section{Proof of Theorem~\ref{1233}}
The proof of  Theorem~\ref{1233} follows immediately from  Theorem~\ref{1612} and Theorem~\ref{1535}. However, the additional assumption (\ref{1427}) significantly simplifies the proof. 

Consider a set of sequences of reals
$$
\{a_{11}\}, \{a_{12}, a_{22}\},\dots, \{a_{1n}, a_{2n},\dots a_{nn}\}, \dots\;,
$$
where all $a_{jn}\in [0,1]$. We say that $\{a_{jn}\}$ are uniformly distributed in $[0,1]$ if for any $0\leqslant a<b\leqslant 1$
$$
\lim_{n\to\infty}\frac{\#\{j\in\{1,2,\cdots ,n\}\,:\,a_{jn}\in [a,b]\}}{n}=b-a\;.
$$ 
The definition is an insignificant generalization of the notion of uniformly distributed sequences (see, e.g., \cite{KN74}). It is easy to see that the Weyl criterion (see Ibid.) continues to be valid in this case:

{\it The set of sequences} $\{a_{jn}, j=1,\dots ,n\},\,n=1,2,\dots,$ {\it is uniformly distributed if and only if for all} $l=1,2,\dots$
$$
\frac{1}{n}\sum_{j=1}^n e^{2\pi la_{jn}}\to 0,\quad n\to\infty\;.
$$

Let $z_{jn}=r_{jn}e^{i\theta_{jn}}$ be a zero of $G_n (z),\,r_{jn}=|z_{jn}|,\,\theta_{jn}=\arg z_{jn},\,0\leqslant\theta_{jn}<2\pi.$ The asymptotic uniform distribution of the arguments is equivalent to the statement that the set of sequences $\{\theta_{jn}/(2\pi)\}$ is uniformly distributed. Thus, according to Weyl's criterion, it is enough to show that for any
$l=1,2,\dots$
$$ \lim_n\frac{1}{n}\sum_{j=1}^n e^{il\theta_{jn}} =0$$
with probability 1.

For the simplicity we assume that $\xi_0\ne0$. Consider the random polynomial
$$
\tilde{G}_n(z)=\xi_n+\xi_{n-1}z+\dots+\xi_{1} z^{n-1} +\xi_0 z^n\;.
$$
Its roots are $z_{kn}^{-1}$. According to Newton's formulas (see, e.g., \cite{aK88}),
$$
\sum_{j=1}^n\frac{1}{z_{jn}^l} =\varphi_l\left (\frac{\xi_1}{\xi_0},\dots \frac{\xi_l}{\xi_0}\right)\;,
$$
where $\varphi_l (x_1,\dots x_l)$ are polynomials which do not depend on $n$. (For example, $\varphi_1(x)=-x$). It follows that
\begin{equation}\label{1540}
\frac{1}{n}\sum_{j=1}^{n}
e^{-il\theta_{jn}}=\frac{1}{n}\sum_{j=1}^n e^{-il\theta_{jn}}\left(1-\frac{1}{r_{jn}^l}\right)+\frac{\varphi_l}{n}\;.
\end{equation}

 As was shown in the proof of Theorem~\ref{1612}, for $|z|<1$ the polynomials $G_n(z)$ converge to the analytical function $G(z)=\sum_{k=0}^{\infty}\xi_k z^k$ with probability 1. Since  $\xi_0\ne0$, the function $G(z)$ has no zeros inside a circle $\{z:|z|\leqslant \rho\},\,
 \mathbf{P}\{\rho >0\}=1$. Hence for $n\geqslant N,\,\mathbf{P}\{N<\infty\},$ the polynomials $G_n(z)$ have no zeros inside $\{z:\, |z|\leqslant\rho\}.$ Let $\gamma >0$ be a positive number. It follows from (\ref{1540}) that
$$
\Big|\frac{1}{n}\sum_{j=1}^n e^{il\theta_{jn}}\Big|\leqslant(l+1)\frac{\gamma}{(1-\gamma )^l}+\frac{1}{n}\left(1+\frac{1}{\rho}\right)\#\{j:|r_{jn}-1|>\gamma, i=1,\dots n\}+\frac{\varphi_l}{n}\;.
$$
Theorem~\ref{1612} implies that the second member on the right-hand side goes to zero as $n\to\infty$ with probability 1. Hence 
$$
\frac{1}{n}\sum_{j=1}^n e^{il\theta{jn}}\to 0,\quad n\to\infty\;,
$$
with probability 1 and the theorem follows.

\section{Proof of Theorem~\ref{1535}}
Consider integer numbers $p,q_1,q_2$ such that $0\leqslant q_1<q_2<p-1$. Put $\varphi_j=q_j/p,\,j=1,2,$ and try to estimate $S_n=S_n(2\pi\varphi_1,2\pi\varphi_2)$. Evidently $S_n=\lim_{R\to\infty}S_{nR}$, where $S_{nR}$ is the number of zeros of $G_n(z)$ inside the domain $A_R=\{z\,:\,|z|\leqslant R, 2\pi\varphi_1\leqslant\arg z\leqslant2\pi\varphi_2\}$. It follows from the argument principle (see, e.g., \cite{aM83}) that $S_{nR}$ is equal to the change of the argument of $G_n(z)$ divided by $2\pi$  as $z$ traverses the boundary of $A_R$. The boundary consists of the arc $\Gamma_R=\{z\,:\,|z|=R,2\pi\varphi_1\leqslant\arg z\leqslant2\pi\varphi_2\}$ and two intervals $L_j=\{z\,:\,0\leqslant|z|\leqslant R,\arg z=\pi\varphi_j\},\,j=1,2$. It can easily be checked that if $R$ is sufficiently large, then the change of the argument as $z$ traverses $\Gamma_R$ is equal to $n(\varphi_2-\varphi_1)+o(1)$ as $n\to\infty$. If $z$ traverses a subinterval of $L_j$ and the change of the argument of $G_n(z)$ is at least $\pi$, then
the function $|G_n(z)|\cos(\arg G_n(z))$ has at least one root in this interval. It follows from Theorem~\ref{1713} that with probability one the number of real roots of the polynomial
$$
g_{n,j}(x)=\sum_{k=0}^nx^k\Re(\xi_ke^{2\pi ik\varphi_j})=\sum_{k=0}^nx^k\eta_{k,j}
$$
is $o(n)$ as $n\to\infty$. Thus the change of the argument of $G_n(z)$ as $z$ traverses $L_j$ is $o(n)$ as $n\to\infty$ and
$$
\P\left\{\frac1nS_n(2\pi\varphi_1,2\pi\varphi_2)=(\varphi_2-\varphi_1)+o(1),\quad n\to\infty\right\}=1\;.
$$

The set of points of the form $\exp\{2\pi iq/p\}$ is dense in the unit circle $\{z\,:\, |z|=1\}$. Therefore 
$$
\P\left\{\frac1nS_n(\alpha,\beta)\underset{n\to\infty}{\longrightarrow}\frac{\beta-\alpha}{2\pi}\right\}=1
$$
for any $\alpha,\beta$ such that 0 $\leqslant\alpha<\beta\leqslant2\pi$.

\section{Proof of Theorem~\ref{1713}}
First we convert the problem of counting of real zeros of $g_n(x)$ to the problem of counting of sign changes in the sequence of the derivatives $\{g_n^{(j)}(1)\}_{j=0}^n$.

Let $\{a_j\}_{j=0}^n$ be a sequence of real numbers. By $Z(\{a_j\})$ denote the number of sign changes in the sequence $\{a_j\}$, which is defined as follows. First we exclude all zero members from the sequence. Then we count the number of the neighboring members of different signs.

For any polynomial $p(x)$ of degree $n$ put $Z_p(x)=Z(\{p^{(j)}(x)\})$, i.e., the number of sign changes in the sequence $p(x), p^\prime(x),\dots,p^{(n)}(x)$.

\begin{lemma}[Budan-Fourier Theorem]\label{1123}
Suppose $p(x)$ is a polynomial such that $p(a),p(b)\ne0$ for some $a<b$. Then the number of the roots of $p(x)$ inside $(a,b)$ does not exceed $Z_p(a)-Z_p(b)$. Moreover, the difference between $Z_p(a)-Z_p(b)$ and the number of the roots is an even number.
\end{lemma}
\begin{proof}
See, e.g., \cite{aK88}.
\end{proof}
\begin{corollary}
 The number of the roots of $p(x)$ inside $[1,\infty)$ does not exceed $Z_p(1)$.
\end{corollary}
\begin{proof}
For all sufficiently large $x$ the sign of $p^{(j)}(x)$ coincides with the sign of the leading coefficient.
\end{proof}

\begin{corollary}
The function  $Z_p(x)$ does not increase.
\end{corollary}

Let us turn back to the random polynomial $g_n(x)$. Here and elsewhere we shall omit the index $n$ when it can be done without ambiguity. By $M_n(a,b)$ denote the number of zeros of $g(x)$ inside the interval $[a,b]$.

First let us prove that 
\begin{equation}\label{2219}
\E Z_g(1)=o(n),\quad n\to\infty\;.
\end{equation}
Fix some $\varepsilon>0$ and $\lambda\in(0,1/2)$. Since the distributions of $\{\eta_j\}$ belong to a finite set, there exists $K=K(\varepsilon)$ such that
\begin{equation}\label{1614}
\sup_{j\in\Z^1}\P\{|\eta_j|\geqslant K\}\leqslant\varepsilon\;.
\end{equation}

Let $I$ be a subset of $\{0,1,\dots,n\}$ consisting of indices $j$ such that $[\lambda n]\leqslant j\leqslant[(1-\lambda)n]$ and $|\eta_j|<K$. Put
$$
g_1(x)=\sum_{j\in I}\eta_jx^j\;,\quad g_2(x)=g(x)-g(x)\;.
$$

Let $\tau_j$ be the indicator of $\{|g_1^{(k)}(1)|\geqslant|g_2^{(k)}(1)|\}$ and $\chi_j$ be the indicator of $\{|\eta_j|\geqslant K\}$. 

\begin{lemma}\label{2143}
Let $a_1,a_1,b_1,b_2$ be real numbers. If $(a_1+a_2)(b_1+b_2)<0$ and $a_2b_2\geqslant0$, then either $|a_1|\geqslant|a_2|$ or $|b_1|\geqslant|b_2|$.
\end{lemma}
\begin{proof}
The proof is trivial. 
\end{proof}

It follows from Lemma~\ref{2143} that
$$
Z_g(1)=Z_{g_1+g_2}(1)\leqslant Z_{g_2}(1)+2\sum_{j=0}^{n}\tau_j\leqslant Z_{g_2}(1)+2\lambda n+2+2\sum_{j=[\lambda n]}^{[(1-\lambda)n]}\tau_j\;.
$$
Owing to the monotonicity of the function $Z_{g_2}(x)$, one has
$$
Z_{g_2}(1)\leqslant Z_{g_2}(0)\leqslant\sum_{j=0}^{n}\chi_j\;.
$$
Hence,
\begin{equation}\label{1757}
Z_g(1)\leqslant 2\lambda n+2+\sum_{j=0}^{n}\chi_j+2\sum_{j=[\lambda n]}^{[(1-\lambda)n]}\tau_j\;.
\end{equation}
Using (\ref{1614}) we have $\E\chi_j=\P\{|\eta_j|\geqslant K\}\leqslant\varepsilon$, therefore,
\begin{equation}\label{2354}
\E Z_g(1)\leqslant 2\lambda n+2+\varepsilon (n+1)+2\E\sum_{j=[\lambda n]}^{[(1-\lambda)n]}\tau_j\;.
\end{equation}

Let us now estimate $\E\tau_j$. Note that $g^{(k)}(x)=\sum_{l=k}^n\eta_lA_{k,l}x^{l-k}$, where $A_{k,l}=l(l-1)\cdots(l-k+1)$. Fix some integer $k$ such that $\lambda n\leqslant k\leqslant(1-\lambda)n$. If $n-1\geqslant j\geqslant k$, then
$$
A_{k,j}\leqslant (1-\lambda)A_{k,j+1}\;,
$$
which implies
$$
A_{k,j}\leqslant A_{k,[(1-\lambda)n]}(1-\lambda)^{[(1-\lambda)n]-j}
$$
for $\lambda n\leqslant k\leqslant j\leqslant(1-\lambda)n$. Consequently,
$$
|g_1^{(k)}(1)|=\Big|\sum_{j\in J,j\geqslant k}\eta_jA_{k,j}\Big|\leqslant KA_{k,[(1-\lambda)n]}\sum_{j=0}^{[(1-\lambda)n]}(1-\lambda)^j\leqslant\frac{K}{\lambda}A_{k,[(1-\lambda)n]}\;.
$$
This yields  that
\begin{multline*}
\E\tau_k=\P\left\{|g_1^{(k)}(1)|\geqslant|g_2^{(k)}(1)|\right\}\\
\leqslant\P\left\{|g_1^{(k)}(1)|\geqslant|g_1^{(k)}(1)+g_2^{(k)}(1)|-|g_1^{(k)}(1)|\right\}\\
=\P\left\{|g^{(k)}(1)|\leqslant2|g_1^{(k)}(1)|\right\}\leqslant\P\left\{|g^{(k)}(1)|\leqslant\frac{2K}{\lambda}A_{k,[(1-\lambda)n]}\right\}\;.
\end{multline*}

For an arbitrary random variable $X$ define the concentration function $Q(h;X)$ as follows:
$$
Q(h;X)=\sup_{a\in\R}\P\{a\leqslant X\leqslant a+h\}\;.
$$
If $X,Y$ are independent random variables, then (see, e.g., \cite{vP95})
$$
Q(h;X+Y)\leqslant\min\left(Q(h;X),Q(h;Y)\right)\;.
$$
Therefore,
\begin{multline}\label{1644}
\E\tau_k\leqslant\P\left\{\frac{|g^{(k)}(1)|}{A_{k,[(1-\lambda)n]}}\leqslant\frac{2K}{\lambda}\right\}\\
\leqslant\P\left\{\frac{g^{(k)}(1)}{A_{k,[(1-\lambda)n]}}\leqslant\frac{2K}{\lambda}\right\}\leqslant Q\left(\frac{2K}{\lambda};\frac{g^{(k)}(1)}{A_{k,[(1-\lambda)n]}}\right)\\
=Q\left(\frac{2K}{\lambda};\sum_{j=k}^n\frac{ A_{k,j}}{A_{k,[(1-\lambda)n]}}\eta_j\right)\leqslant Q\left(\frac{2K}{\lambda};\sum_{j=[(1-\lambda)n]}^n\frac{ A_{k,j}}{A_{k,[(1-\lambda)n]}}\eta_j\right)\;.
\end{multline}
To estimate the right-hand side of (\ref{1644}) we use the following result.
\begin{lemma}[the Kolmogorov-Rogozin inequality]
Let $X_1,X_2,\dots,X_n$ be independent random variables. Then for any $0<h_j\leqslant h,\,j=1,\dots,n,$
\begin{equation}\label{1714}
Q(h;X_1+\cdots+X_n)\leqslant\frac{Ch}{\sqrt{\sum_{j=1}^nh_j^2(1-Q(h_j;X_j))}}\;,
\end{equation}
where $C$ is an absolute constant.
\end{lemma}
\begin{proof}
See \cite{bR61}.
\end{proof}

Since the distributions of $\{\eta_j\}$ belong to a finite set, we get
$$
\delta=\delta(\varepsilon,\lambda)=\inf_{j\in\Z^1}\left\{1-Q\left(\frac{2K}{\lambda};\eta_j\right)\right\}>0\;.
$$
Putting $h=h_j=2K/\lambda$ in (\ref{1714}) and using (\ref{1644}), we obtain
\begin{multline*}
\E\tau_k\leqslant C\left[\sum_{j=[(1-\lambda)n]}^n\left\{1-Q\left(\frac{2K}{\lambda};\frac{ A_{k,j}}{A_{k,[(1-\lambda)n]}}\eta_j\right)\right\}\right]^{-1/2}\\
\leqslant C\left[\sum_{j=[(1-\lambda)n]}^n\left\{1-Q\left(\frac{2K}{\lambda};\eta_j\right)\right\}\right]^{-1/2}\leqslant\frac{C}{\sqrt{\delta\lambda n}}\;.
\end{multline*}
Combining this with (\ref{2354}), we have
$$
\E Z_g(1)\leqslant 2\lambda n+2+\varepsilon (n+1)+\frac{2C}{\sqrt{\delta(\varepsilon,\lambda)\lambda}}n^{1/2}\;.
$$
Since $\lambda,\varepsilon$ are arbitrary positive numbers, we obtain (\ref{2219}), which together with the corollary from Lemma~\ref{1123} implies
$$
\E M_n(1,\infty)=o(n)\;,\quad n\to\infty\;.
$$
Considering the random polynomials $g(1/x)$ and $g(-x)$, it is possible to obtain similar estimates for $M_n(0,1)$ and $M_n(-\infty,0)$. Thus the second part of (\ref{1611}) holds. To prove the first one, we estimate the probabilities of large deviations for the sums $\sum\chi_j$ and $\sum\tau_j$. The elementary considerations or the application of Bernstein inequalities (see, e.g., \cite{vP95}) leads to
\begin{equation}\label{1755}
\P\left\{\Big|\sum_{j=0}^n\chi_j\Big|>2(n+1)\varepsilon\right\}\leqslant2e^{-n\varepsilon/8}\;.
\end{equation}
The analysis of the behavior of $\sum\tau_j$ is slightly more difficult . 

Henceforth we shall use the following notation: for any positive functions $f_1,f_2$ we write $f_1\ll f_2$, if there exists an absolute constant $C$ such that $f_1\leqslant Cf_2$ in the domain of these functions.

\begin{lemma}\label{1228}
There exists a constant $c$ depending only on $\lambda,\varepsilon$ and the distributions of $\{\eta_j\}$ such that
$$
\E\tau_k\leqslant cn^{-2}
$$
for $\lambda n\leqslant k\leqslant(1-\lambda)n$.
\end{lemma}
\begin{proof}
As was shown in (\ref{1644}),
\begin{equation}\label{1833}
\E\tau_k\leqslant Q\left(\frac{2K}{\lambda};\sum_{j=[(1-\lambda)n]}^n\frac{ A_{k,j}}{A_{k,[(1-\lambda)n]}}\eta_j\right)\;\;.
\end{equation}
\end{proof}
To estimate the concentration function in the right-hand side we use the result of Esseen (see, e.g., \cite{vP95}). Let $X$ be a random variable with a characteristic function $f(t)$. Then 
$$
Q(h;X)\ll\max\left(h,\frac1T\right)\int_{-T}^T|f(t)|\,dt
$$
uniformly for all $T>0$.

Putting $T=\lambda/(KA_{k,[(1-\lambda)n]})$ and applying (\ref{1833}) , we obtain
$$
\E\tau_k\ll\frac 1T\int_{-T}^T\prod_{j=[(1-\lambda)n]}^n|f_j(A_{kj}t)|\,dt\;,
$$
where $f_j(t)$ is a characteristic function of $\eta_j$. Further,
\begin{multline*}
\E\tau_k\ll\frac 1T\int_{-T}^T\left[\prod_{j=[(1-\lambda)n]}^n|f_j(A_{kj}t)|^2\right]^\frac12\,dt\\
\ll\frac 1T\int_{-T}^T\exp\left\{-\frac12\sum_{j=[(1-\lambda)n]}^n\left(1-|f(A_{kj}t)|^2\right)\right\}\,dt
\\=\frac 1T\int_{-T}^T\exp\left\{-\frac12\sum_{j=[(1-\lambda)n]}^n\int_{-\infty}^\infty\left[1-\cos(A_{kj}tx)\right]\,\mathcal{P}_j(dx)\right\}\,dt\;,
\end{multline*}
where $\mathcal{P}_j$ is a distribution of the symmetrized $\eta_j$, i.e., a distribution of $\eta_j-\eta^\prime_j$, where $\eta^\prime_j$ is an independent copy of $\eta_j$.

There are at most $r$ different distributions among $\{\mathcal{P}_j\}_{(1-\lambda)n\leqslant j\leqslant n}$. Therefore there exist a distribution $\mathcal{P}$ and a subset $J\subset\{j\,:\,(1-\lambda)n\leqslant j\leqslant n\}$ such that $|J|\geqslant n\lambda/r$ and $\mathcal{P}_j=\mathcal{P}$ for all $j\in J$. By $\sum^\prime$ denote the summation taking over all indices such that $j\in J$. Thus,
$$
\E\tau_k\ll\frac 1T\int_{-T}^T\exp\left\{-\frac12\sideset{}{'}\sum_{j=[(1-\lambda)n]}^n\int_{-\infty}^\infty\left[1-\cos(A_{kj}tx)\right]\,\mathcal{P}(dx)\right\}\,dt\;.
$$
Choose $\delta>0$ such that $\gamma=\mathcal{P}\{x\,:\,|x|>\delta\}>0$. Since the integrands are non-negative, we get
\begin{multline*}
\E\tau_k\ll\frac 1T\int_{-T}^T\exp\left\{-\frac12\sideset{}{'}\sum_{j=[(1-\lambda_r)n]}^n\int_{|x|>\delta}\left[1-\cos(A_{kj}tx)\right]\,\mathcal{P}(dx)\right\}\\
=\frac 1T\int_{-T}^Te^{-\beta n+s(t)}\,dt\;,
\end{multline*}
where $\lambda_r=\lambda(2r-1)/(2r),\,\beta=|J\cap\{j\,:\,(1-\lambda_r)n\leqslant j\leqslant n\}|/(2n)$ and 
$$
s(t)=\frac12\int_{|x|>\delta}\sideset{}{'}\sum_{j=[(1-\lambda_r)n]}^n\cos(A_{kj}tx)\,\mathcal{P}(dx)\;.
$$
Put $\alpha=\lambda\gamma/(4r)$ and consider $\Lambda_1=\{t\in[-T,T]\,:\,|s(t)|<\alpha n/2\}$ and $\Lambda_2=[-T,T]\setminus\Lambda_1$. Since $|J|\geqslant n\lambda/r$ and by the definition of $\beta$, we have $\beta\geqslant\alpha$. Therefore,
\begin{equation}\label{2205}
\E\tau_k\ll e^{-\alpha n/2}+\frac{\mu(\Lambda_2)}{T}\;,
\end{equation}
where  $\mu$ denotes the Lebesgue measure.

Let us estimate $\mu(\Lambda_2)$. It follows from Chebyshev's and H\"older's inequalities that
\begin{equation}\label{2209}
\mu(\Lambda_2)\leqslant\frac{16}{\alpha^4n^4}\int_{-T}^T|s(t)|^4\,dt\leqslant\frac{1}{\alpha^4n^4}\int_{|x|>\delta}\,d\mathcal{P}\int_{-T}^T\Big|\sideset{}{'}\sum_{j=[(1-\lambda_r)n]}^{n}\cos(A_{kj}tx)\Big|^4\,dt\;.
\end{equation}
Put
$$
S(x)=\int_{-T}^T\Big|\sideset{}{'}\sum_{j=[(1-\lambda_r)n]}^{n}\cos(A_{kj}tx)\Big|^4\,dt
$$
and assume, for simplicity, that $r=1$, i.e., $\lambda_r=\lambda/2,\,\sum=\sum^\prime$ and the summation is taken over all $j$.  The general case is considered in a similar way.

We have
\begin{multline}\label{22258}
S(x)=\int_{-T}^T\bigg(\sum_{j_1}\cos^4(A_{kj_1}tx)+\sum_{j_1\ne j_2}\cos^3(A_{kj_1}tx)\cos(A_{kj_2}tx)\\
+\sum_{j_1\ne j_2}\cos^2(A_{kj_1}tx)\cos^2(A_{kj_2}tx)\\
+\sum_{j_1\ne j_2\ne j_3}\cos^2(A_{kj_1}tx)\cos(A_{kj_2}tx)\cos(A_{kj_3}tx)\\
+\sum_{j_1\ne j_2\ne j_3\ne j_4}\cos(A_{kj_1}tx)\cos(A_{kj_2}tx)\cos(A_{kj_3}tx)\cos(A_{kj_4}tx)\bigg)\,dt\;.
\end{multline}
The first three summands in (\ref{22258})  are easily estimated as follows:
\begin{multline}\label{1841}
\Big|\int_{-T}^T\bigg(\sum_{j_1}\cos^4(A_{kj_1}tx)+\sum_{j_1\ne j_2}\cos^3(A_{kj_1}tx)\cos(A_{kj_2}tx)\\
+\sum_{j_1\ne j_2}\cos^2(A_{kj_1}tx)\cos^2(A_{kj_2}tx)\bigg)\,dt\Big|\ll Tn^2\;.
\end{multline}
The next two summands have a common method of estimation. We consider only the last one. From the formula $\cos y=(e^{iy}+e^{-iy})/2$ it is easily shown that
\begin{multline}\label{1217}
\Big|\int_{-T}^T\sum_{j_1\ne j_2\ne j_3\ne j_4}\cos(A_{kj_1}tx)\cos(A_{kj_2}tx)\cos(A_{kj_3}tx)\cos(A_{kj_4}tx)\,dt\Big|\\
\ll\sum_{j_1\ne j_2\ne j_3\ne j_4}\min\left(T,|x|^{-1}|\pm A_{kj_1}\pm A_{kj_2}\pm A_{kj_3}\pm A_{kj_4}|^{-1}\right)\\
\ll\sum_{j_1>j_2>j_3>j_4}\min\left(T,|x|^{-1}A_{kj_1}^{-1}\Big|1-\frac{A_{kj_2}}{A_{kj_1}}-\frac{A_{kj_3}}{A_{kj_1}}-\frac{A_{kj_4}}{A_{kj_1}}\Big|^{-1}\right)\;.
\end{multline}
The summation in the middle term is taken over all possible combinations of signs. 

Consider the partition of the index set $\{j=(j_1,j_2,j_3,j_4)\,:\,j_1>j_2>j_3>j_4\}=K_1\cup K_2$, where
$$
K_1=\left\{j\,:\,j_1-j_2\leqslant\frac{10}{\lambda},\, j_1-j_3\leqslant\frac{10}{\lambda}|\ln\lambda|\right\}
$$
and $K_2$ is the complement of $K_1$. Clearly, $|K_1|\ll n^2|\ln\lambda|/\lambda^2$. Therefore,
\begin{equation}\label{1218}
\sum_{j\in K_1}\min\left(T, |x|^{-1}A_{kj_1}^{-1}\Big|1-\frac{A_{kj_2}}{A_{kj_1}}-\frac{A_{kj_3}}{A_{kj_1}}-\frac{A_{kj_4}}{A_{kj_1}}\Big|^{-1}\right)\ll\frac{Tn^2|\ln\lambda|}{\lambda^2}\;.
\end{equation}
Consider now
$$
\sum_{j\in K_2}A_{kj_1}^{-1}\Big|1-\frac{A_{kj_2}}{A_{kj_1}}-\frac{A_{kj_3}}{A_{kj_1}}-\frac{A_{kj_4}}{A_{kj_1}}\Big|^{-1}\;.
$$
Putting $p=j_1-j_2$, we have
\begin{multline*}
\frac{A_{kj_2}}{A_{kj_1}}=\frac{(j_1-p)\cdots(j_1-p-k+1)}{j_1\cdots(j_1-k+1)}\\
=\left(1-\frac{p}{j_1}\right)\cdots\left(1-\frac{p}{j_1-k+1}\right)\leqslant\exp\left\{-p\sum_{l=j_1-k+1}^{j_1}\frac1l\right\}.
\end{multline*}
Since for any natural $l$
$$
\frac{1}{l}>\ln\left(1+\frac{1}{l}\right)=\ln(l+1)-\ln l\;,
$$
we get
$$
\sum_{l=j_1-k+1}^{j_1}\frac1l>\ln(j_1+1)-\ln(j_1-k+1)=-\ln\left(1-\frac{k}{j_1+1}\right)\;.
$$
Taking into account $\lambda n\leqslant k\leqslant(1-\lambda)n$ and $ (1-\lambda/2)n\leqslant j_1\leqslant n$ and using the inequality 
$$
-\ln(1-t)\geqslant t\;,\quad t\in[0,1]\;,
$$
we get
$$
\sum_{l=j_1-k+1}^{j_1}\frac1l\geqslant\frac{\lambda n}{n+1}\geqslant\frac12\lambda\;.
$$
Therefore,
\begin{equation}\label{1640}
\frac{A_{kj_2}}{A_{kj_1}}\leqslant\exp\left\{-\frac \lambda2p\right\}=\exp\left\{-\frac\lambda2(j_1-j_2)\right\}\;.
\end{equation}

If $j\in K_2$ and $j_1-j_2>10/\lambda$, then 
$$
\frac{A_{kj_4}}{A_{kj_1}}\leqslant\frac{A_{kj_3}}{A_{kj_1}}\leqslant\frac{A_{kj_2}}{A_{kj_1}}\leqslant e^{-5}<\frac14\;,
$$
which implies
\begin{equation}\label{1650}
1-\frac{A_{kj_2}}{A_{kj_1}}-\frac{A_{kj_3}}{A_{kj_1}}-\frac{A_{kj_4}}{A_{kj_1}}\geqslant\frac14\;.
\end{equation}

Suppose now $j\in K_2$ and $j_1-j_3>10|\ln\lambda|/\lambda$. Using (\ref{1640}) and $\lambda\in(0,1/2)$, we get 
$$
1-\frac{A_{kj_2}}{A_{kj_1}}\geqslant1-e^{-\lambda/2}\geqslant\frac\lambda2\left(1-\frac\lambda4\right)\geqslant\frac{7}{16}\lambda\;.
$$
Further, (\ref{1640}) also holds for $j_3$. Therefore,
$$
\frac{A_{kj_4}}{A_{kj_1}}\leqslant\frac{A_{kj_3}}{A_{kj_1}}\leqslant\exp\left\{-\frac\lambda2(j_1-j_3)\right\}\leqslant\exp\left\{-\frac{10}{2}|\ln\lambda|\right\}\leqslant\lambda^5\leqslant\frac{1}{16}\lambda\;.
$$
Thus,
\begin{equation}\label{1707}
1-\frac{A_{kj_2}}{A_{kj_1}}-\frac{A_{kj_3}}{A_{kj_1}}-\frac{A_{kj_4}}{A_{kj_1}}\geqslant\frac{5}{16}\lambda\;.
\end{equation}
It follows from (\ref{1650}) and (\ref{1707}) that 
$$
\sum_{j\in K_2}A_{kj_1}^{-1}\Big|1-\frac{A_{kj_2}}{A_{kj_1}}-\frac{A_{kj_3}}{A_{kj_1}}-\frac{A_{kj_4}}{A_{kj_1}}\Big|^{-1}\ll\frac{1}{\lambda}\sum_{j}A_{kj_1}^{-1}\;.
$$
Taking into account the structure of the index set $\{j\}$, we have
$$
\sum_{j}A_{kj_1}^{-1}\leqslant \frac{(\lambda n)^4}{A_{k,[(1-\lambda/2)n]}}\;,
$$
consequently,
\begin{equation}\label{1721}
\sum_{j\in K_2}A_{kj_1}^{-1}\Big|1-\frac{A_{kj_2}}{A_{kj_1}}-\frac{A_{kj_3}}{A_{kj_1}}-\frac{A_{kj_4}}{A_{kj_1}}\Big|^{-1}\ll\frac{\lambda^3n^4}{A_{k,[(1-\lambda/2)n]}}\;.
\end{equation}

Combining (\ref{22258}), (\ref{1841}), (\ref{1217}), (\ref{1218}) and (\ref{1721}), we obtain
$$
S(x)\ll Tn^2+\frac{Tn^2|\ln\lambda|}{\lambda^2}+\frac{\lambda^3n^4}{|x|A_{k,[(1-\lambda/2)n]}}\;.
$$
Applying this to (\ref{2209}), we get
$$
\mu(\Lambda_2)\ll \frac{T}{\alpha^{4}n^{2}}+\frac{T|\ln\lambda|}{\lambda^2\alpha^{4}n^{2}}+\frac{\lambda^3}{\alpha^4\delta A_{k,[(1-\lambda/2)n]}}\;.
$$
By (\ref{2205}),
$$
\E\tau_k\ll e^{-\alpha n/2}+\frac{1}{\alpha^{4}n^{2}}+\frac{|\ln\lambda|}{\lambda^2\alpha^{4}n^{2}}+\frac{\lambda^3}{T\alpha^4\delta A_{k,[(1-\lambda/2)n]}}\;.
$$
Recalling that  $T=\lambda/(KA_{k,[(1-\lambda)n]})$, we obtain
$$
\E\tau_k\ll e^{-\alpha n/2}+\frac{1}{\alpha^{4}n^{2}}+\frac{|\ln\lambda|}{\lambda^2\alpha^{4}n^{2}}+\frac{\lambda^2KA_{k,[(1-\lambda)n]}}{\alpha^4\delta A_{k,[(1-\lambda/2)n]}}\;.
$$
It follows from (\ref{1640}) that
$$
\frac{A_{k,[(1-\lambda)n]}}{A_{k,[(1-\lambda/2)n]}}\leqslant e^{-\lambda^2n/4}\;.
$$
Thus,
$$
\E\tau_k\ll e^{-\alpha n/2}+\frac{1}{\alpha^{4}n^{2}}+\frac{|\ln\lambda|}{\lambda^2\alpha^{4}n^{2}}+\frac{\lambda^2K}{\alpha^4\delta }e^{-\lambda^2n/4}\;.
$$
Recalling that  $\alpha=\gamma\lambda/4$, we obtain
$$
\E\tau_k\ll e^{-\gamma\lambda n/8}+\frac{1}{\gamma^{4}\lambda^4n^{2}}+\frac{|\ln\lambda|}{\gamma^4\lambda^6n^{2}}+\frac{K}{\gamma^4\lambda^2\delta }e^{-\lambda^2n/4}\;.
$$
Since $K$ is defined by $\varepsilon$ and $\gamma,\delta$ are defined by the distributions of $\{\eta_j\}$, Lemma~\ref{1228} is proved.

Now we are ready to complete the proof of Theorem~\ref{1713}. It follows from  (\ref{1757}) that
\begin{equation}\label{1313}
M_n(1,\infty)\leqslant 2\lambda n+2+\sum_{j=0}^{n}\chi_j+2\sum_{j=[\lambda n]}^{[(1-\lambda)n]}\tau_j\;.
\end{equation}
By  Lemma~\ref{1228} and  Chebyshev's inequality,
\begin{equation}\label{1314}
\P\left\{\sum_{k=[\lambda n]}^{[(1-\lambda)n]}\tau_k>n^{3/4}\right\}\leqslant\frac{\sum_{j=[\lambda n]}^{[(1-\lambda)n]}\E\tau_k}{n^{3/4}}\leqslant c_1n^{-5/4}\;.
\end{equation}
Further, it follows from (\ref{1755}) that there exists a constant $c_2>0$ depending only on $\varepsilon$ such that
\begin{equation}\label{1315}
\P\left\{\sum_{j=0}^n\chi_j>2\varepsilon n\right\}\leqslant c_2n^{-2}\;.
\end{equation}
Combining (\ref{1313}), (\ref{1314}) and (\ref{1315}), we get
$$
\P\left\{M_n(1,\infty)>2\lambda n+2+2n^{3/4}+2\varepsilon n\right\}\leqslant c_1n^{-5/4}+c_2n^{-2}\;.
$$
Considering the random polynomials $g(1/x)$ and $g(-x)$, it is possible to obtain similar estimates for $M_n(0,1)$ and $M_n(-\infty,0)$. Thus there exist positive constants $c'_1,c'_2$ such that
$$
\P\left\{M_n>2\lambda n+2+2n^{3/4}+2\varepsilon n\right\}\leqslant c'_1n^{-5/4}+c'_2n^{-2}\;.
$$

According to the Borel-Cantelli  lemma, with probability one there exists only a finite number of $n$ such that $M_n>2\lambda n+2+2n^{3/4}+2\varepsilon n$. Since $\lambda,\varepsilon$ are arbitrary small,
$$
\P\left\{\frac{M_n}{n}\underset{n\to\infty}{\longrightarrow}0\right\}=1\;.
$$
Theorem~\ref{1713} is proved.

\emph{Acknowledgements} A part of the work was done in the University of Bielefeld. The authors thank F.~G\"otze for the possibility to participate at the work of CRC 701 ``Spectral Structures and Topological Methods in Mathematics''. They are also grateful to A.~Cole  for her valuable help.

\end{document}